\documentclass[journal=iecred]{achemso}

\usepackage[version=4]{mhchem} 
\usepackage{amsmath}
\usepackage{amssymb}
\usepackage{siunitx}
\SectionNumbersOn

\author{Robert J. Lovelett}
\affiliation{Department of Chemical and Biological Engineering, Princeton University}
\author{Jos\'{e} L. Avalos}
\affiliation{Department of Chemical and Biological Engineering, Princeton University}
\author{Ioannis G. Kevrekidis}
\affiliation{Department of Chemical and Biological Engineering, Princeton University}
\alsoaffiliation{Department of Chemical and Biomolecular Engineering, Johns Hopkins University}
\email{yannis@princeton.edu}
\title[Grey-Box Artificial Neural networks]
  {Partial observations and conservation laws: Grey-box modeling in biotechnology and optogenetics}

\begin{document}

\begin{tocentry}




    \centering
    \includegraphics{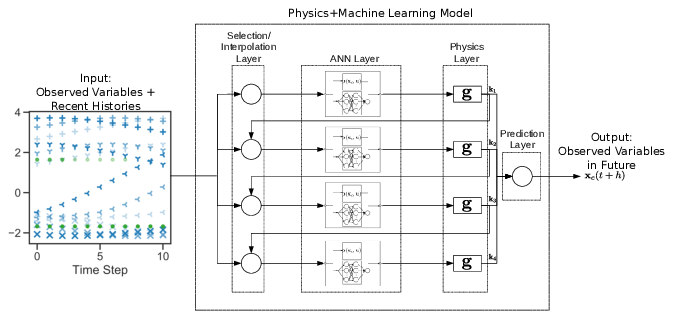}

\end{tocentry}

\begin{abstract}
Developing accurate dynamical system models from physical insight or data can be impeded when only partial observations of the system state are available.
Here, we combine conservation laws used in physics and engineering with artificial neural networks to construct ``grey-box'' system models that make accurate predictions even with limited information.
These models use a time delay embedding (c.f., Takens embedding theorem) to reconstruct effect of the intrinsic states, and can be used for multiscale systems where macroscopic balance equations depend on unmeasured micro/meso scale phenomena.
By incorporating physics knowledge into the neural network architecture, we regularize variables and may train the model more accurately on smaller data sets than black-box neural network models.
We present numerical examples from biotechnology, including a continuous bioreactor actuated using light through optogenetics (an emerging technology in synthetic biology) where the effect of unmeasured intracellular information is recovered from the histories of the measured macroscopic variables.

\end{abstract}

\section{Introduction}
Developing a system model is the necessary first task for effective process control or optimization.
In general, one can develop a system model using one of two philosophies: a theory-based, physics-driven model, or an empirically structured, data-driven model.
The latter approach (``system identification'') can be further subdivided into a purely data-driven black-box approach, or a physics informed grey-box approach \cite{Ljung}.
In this work, we examine a nonlinear grey-box approach for system identification that uses both artificial neural networks and physics-based modeling equations.
The known information about the underlying system (the ``greyness'') typically arises from application-specific macroscopic conservation laws that are widely used in physics and engineering, while the unknown information typically involves constitutive laws (``closures'') that relate the effects of the micro-/meso- scale state of the system to the macro-scale variables, \textit{e.g.}, the effect of intracellular concentrations on the apparent microbial kinetics.
We focus in particular on the example of a bioreactor actuated by light inputs that control microbial gene expression using an emerging technology called optogenetics\cite{Toettcher2011,Milias-argeitis2015,Milias-Argeitis2016,Zhao2018,Lalwani2018}.
Our grey-box approach is well-suited for this problem given that it is straightforward to construct a set of physics-based balance equations that constrain the dynamics, yet full mechanistic models that include all relevant microscale/intracellular variables are challenging to develop and expensive to validate. 

Artificial neural networks (ANNs) have been widely used by the control systems community owing partly to their property that they can approximate any continuous function in Euclidean space\cite{Narendra1990}.
ANNs are especially useful in modeling nonlinear dynamical systems, where features such as bifurcation and saturation are important to consider for designing the control system; this has led to their use in black-box system identification and control for decades \cite{Narendra1990, Lee2018, Venkatasubramanian2019}.
ANNs have re-emerged more recently in the literature partly due to of a number of high-profile and promising results in other domains, such as models for classifying images\cite{Krizhevsky2012} or for ``playing'' complex games like Go more effectively than expert human players\cite{Silver2017}.
These advances have been enabled by the use of ``deep'' nets (neural networks with many hidden layers), which are now possible to train effectively because of increased computational power and better stochastic training algorithms that can handle terabyte-scale data sets\cite{Lecun2015}.
Despite these advances, several authors have recommended using caution when applying deep neural nets (and other modern machine learning tools) for process systems engineering applications---including raising concerns that performance will not be improved in comparison with traditional tools \cite{Rawlings2019} and that they lack the rigor, explicability, and generalization properties of either simpler black-box or first-principles models.\cite{Venkatasubramanian2019}
These issues motivate developing methods that can build on the recent advances in deep neural nets, while preserving certain properties of more traditional modeling tools. 

To alleviate some of these concerns, we investigate a hybrid modeling approach that includes both first-principles knowledge and ``deep'' data-driven neural network models, often referred to as grey-box ANNs (GB-ANNs).
The neural network models in the GB-ANN are used as ``constitutive laws'' (e.g., chemical reaction rate laws) within known, physics-based conservation laws that govern the system dynamics.
By including knowledge of conservation law principles into the architecture of the GB-ANNs, the model outputs are strongly regularized in comparison with black-box ANN (BB-ANN) models, which enables training on more modest sized data sets such as those collected in scientific research or pilot-scale industrial systems.

The general philosophy of combining artificial neural networks with structured physics-based equations has been investigated previously by several authors.
Psichogios and Ungar developed a discrete-time hybrid model of a batch bioreactor that used a neural network model for cell growth, and they observed a reduction in prediction error when they compared their hybrid model with a black-box neural network model.\cite{Psichogios1992}
Rico-Martinez, Anderson, and Kevrekidis used a continuous time framework to demonstrate how a grey-box neural network model of a dynamical system learned from data can be used for bifurcation analysis, which may be difficult in the discrete time scenario.\cite{Rico-Martinez1994}
More recently, Raissi, Perdikaris, and Karniadakis used what they refer to as ``physics-informed neural networks'' to learn unknown state-dependent parameters for a class of systems governed by nonlinear partial differential equations that arise in fluid flow, reaction-diffusion, and other problems; they observe that, by incorporating physical knowledge, they can regularize the outputs and closely predict out-of-sample data.\cite{Raissi2019}

One significant weakness of the GB-ANN models presented in the literature to date is the requirement that the unknown component of the dynamics (which we call the constitutive law) is a function strictly of the modeled variables, when in general it may depend on unknown or unmeasured, often micro-scale, state variables as well.
In the context of black-box modeling, this weakness is often overcome by reconstructing the state space using a delay embedding of the observed and/or modeled variables; this approach was first introduced in the dynamical systems field by Takens,\cite{Takens1981} relying on the Whitney embedding theorem \cite{Whitney}.
The delay embedding approach has been used to find accurate continuous-time nonlinear models of complex spatially distributed reacting systems from experimental data \cite{Rico-Martinez1992,Krischer1993}.
More recently, parametrizing systems using histories of observations has been used by several authors to ``learn'' physical laws from data in combination with several machine learning algorithms\cite{Bongard2007, Schmidt2009, Brunton2015, Yair2017}.
Here, we will parametrize the state of the system using the known macroscopic physical variables and their recent histories, and the final system model will be a delay differential equation that can be written in terms of only known macroscopic variables.

Often, phenomenological approximations of constitutive laws (e.g., enzyme kinetics laws, such as the Michaelis-Menten equation\cite{Johnson2011}) can be an informative starting point for developing a dynamical model (in fact, they are sometimes sufficient on their own).
When these models are available (with or without parameter estimates), we can augment the GB-ANN to use this information.
The GB-ANN in this case will \emph{learn the difference} between the phenomenological estimate and the underlying true constitutive relation.
By inserting these relationships into the ANN model, we can potentially reduce the amount of data required to train an accurate model, allowing for more rapid process development.

One particularly promising area where we believe our approach will find use is in designing control systems for bioreactors used for microbial production of fuels or chemicals that are actuated by light using optogenetics\cite{Toettcher2011,Milias-argeitis2015,Milias-Argeitis2016,Lalwani2018, Zhao2018}.
These systems typically work by using a light-activated transcription factor to express or repress a set of genes, such as those that code for enzymes that catalyze production of a biofuel\cite{Zhao2018}.
Several research groups have proposed and are developing \textit{in silico} control systems for light-controlled microbial systems\cite{Milias-argeitis2015,Milias-Argeitis2016,Lalwani2018}, but there has been relatively little development of the accurate dynamical models required for optimization and control, whether from data or from first principles.
Now, while it is straightforward to construct a set of macroscopic balance equations that describe the time evolution of biomass, nutrients, and products, \cite{Bailey} the rate laws that govern these equations are unlikely to be available.
Furthermore, typical phenomenological assumptions, such as the Monod equation for cell growth rate\cite{Monod1949}, are unlikely to be valid when the microbes are driven to different microscopic/intracellular states by varying the light input.
Deep ANN models, however, can be trained using time delay reconstruction of the partially unmeasured system state to serve as effective rate laws that can be learned from a modest amount of data using only macroscopic variables and their histories.

The remainder of this paper is organized as follows: in Section \ref{sec:theory}, we present a motivating example, the problem formulation, GB-ANN structure, and the system identification (i.e., training) methodology; then, in Section \ref{sec:res_disc}, we discuss two illustrative numerical examples---(1) as a validation step, a simple demonstration using a model of chemostat that implements a biochemical reaction, and (2) a more realistic example of a model of a continuous bioreactor actuated by light via optogenetics; finally in Section \ref{sec:conc}, we present conclusions and suggestions for future research.

\section{Theory}\label{sec:theory}
Before a formal mathematical treatment, we present as a motivating example the case of a continuous stirred bioreactor (i.e., a chemostat) where a microbial species is growing, consuming nutrients, and secreting products.
It is straightforward to write a macroscopic conservation law on the bioreactor that describes the dynamics of biomass (scalar $X$) and extracellular chemical species (vector $\mathbf{C}$), for example:
\begin{equation} \label{eqn:simple_chemo_balance}
    \begin{aligned}
        \dot X &= (\mu(X,\mathbf{C}) - D)X\\
        \dot C_i &= k_i(X,\mathbf{C}) X + (C_{i0} - C_i)D\\
    \end{aligned}
\end{equation}
where $X$ is the concentration of biomass, $C_i$ is the extracellular concentration of chemical species $i$ (where $i$ indexes nutrients and secreted products), $C_{i0}$ is the concentration of chemical $i$ in the feed stream, $\mu$ is the specific growth rate, $D$ is the dilution rate (or the volume-normalized feed and outlet flow rate), and $k_i$ is the secretion (or uptake) rate of species $i$, and  we have assumed that there are no extracellular chemical reactions, no biomass in the feed stream, and density is constant.
In general, the quantities $\mu$ and $k_i$ are unknown functions and are often approximated using phenomenological expressions like the Monod equation for growth rate, $\mu=\frac{\mu_m S}{K_s+S}$ (where parameter $\mu_m$ is the maximum growth rate and parameter $K_s$ is the half saturation constant, and $S$, which is an element of $\mathbf{C}$, is the limiting nutrient/substrate) \cite{Monod1949}.
Such simplified expressions, however, are limited because they only account for the observed extracellular variables ($X$ and $\mathbf{C}$) on which we wrote the balance equations, whereas the actual growth/secretion/uptake rates may depend on the intracellular state of the microbes.

Suppose that there are one or more actuators (``knobs'') that the operator has the ability to manipulate in time, and that these actuators influence the internal metabolism of the microbes.
Here, we consider that the microbe has been engineered using optogenetics to respond to light and will use light intensity as an actuator; with changing light intensity, gene expression in the microbe is altered, resulting in changes the growth rate and to the production (secretion) and uptake rates of chemicals.
Now that the intracellular state of the microbes (the concentrations of proteins, metabolites, etc.) is changing in time  and responding to light, it is unreasonable to assume that  $\mu$ and $k_i$ are functions only of the extracellular variables that are measured.
In this paper, \emph{we will develop a methodology for reconstructing the intracellular state of the system using only extracellular observations} by using the recent history of the observed quantities to infer (``observe'') the effect of the necessary unmeasured/unobserved quantities in Equation \ref{eqn:simple_chemo_balance}.
With these ideas in mind, we will formally and generically define the problem, and then return to this example in Section \ref{sec:cont_bio}.

\subsection{Problem Formulation}
Consider the following state space model of a nonlinear system:
\begin{equation}\label{eqn:sys_general}
\dot{\mathbf{x}} = \mathbf{f}(\mathbf{x},\mathbf{u})
\end{equation}
where $\mathbf{f} \colon \mathbb{R}^{N+M}\to\mathbb{R}^N$ is a smooth, continuous function describing the dynamics, $\mathbf{x}\in \mathbb{R}^N$ is the state vector, and $\mathbf{u}\in \mathbb{R}^M$ is the input actuator vector.
In Equation \ref{eqn:sys_general}, the function $\mathbf{f}$, state $\mathbf{x}$, and their dimensionality $N$ are all unknown.
Now, suppose that a subset of the state variables can be measured (though not necessarily in real-time) and are known quantities that participate in conservation laws on which a set of balance equations may be written, analogous to Equation \ref{eqn:simple_chemo_balance}:
\begin{equation} \label{eqn:sys_conserved}
    \dot{\mathbf{x}}_c = \mathbf{g}(\mathbf{x}_c,\mathbf{u},\pmb{\theta}(\mathbf{x},\mathbf{u})|\pmb{\alpha},\pmb{\beta})
\end{equation}
where $\mathbf{x}_c\in\mathbb{R}^{N_c}$ is a vector containing the set of state variables that participate in conservation laws ($N_c \le N$), $\mathbf{g} \colon \mathbb{R}^{2N_c+M}\to\mathbb{R}^{N_c}$ is a \emph{known} function with \emph{known} parameters $\pmb{\alpha}$ and \emph{unknown} parameters $\pmb{\beta}$, and $\pmb{\theta} \colon \mathbb{R}^{N+M}\to\mathbb{R}^{N_c}$ is an \emph{unknown} function of the \emph{entire} (partially unknown/unmeasured) state vector.
We will refer to $\pmb{\theta}$ as a ``constitutive law'' because it is a system-specific function that is required to write the dynamics of $\mathbf{x}_c$.

Without direct access to the state variables $\mathbf{x}$, we introduce the augmented observable $\mathbf{z}=\left[\mathbf{u}, \mathbf{x}_c\right]^T$ and the following delay embedding of $\mathbf{z}$:
\begin{equation}\label{eqn:delay_embedding_definition}
    \mathbf{z}_e(t)=\begin{bmatrix} \mathbf{z}(t),& \mathbf{z}(t-\tau), & \mathbf{z}(t-2\tau), & ... &  \mathbf{z}(t-(d-1)\tau)\end{bmatrix}^T
\end{equation}
where $\tau\in\mathbb{R}^+$ is the delay time and $d\in\mathbb{N}^+$ is the embedding dimension (see Kim, Eykhalt, and Salas\cite{kim1999} for a discussion of how to select an appropriate delay time).
Subject to weak genericity requirements (see Takens\cite{Takens1981} and Stark et al.\cite{Stark1997, Stark1999, Stark2003}), there exists a diffeomorphism between $\mathbf{x}(t)$ and $\mathbf{z}_e(t)$ as long as $d\ge (2N+1)$ (note that this expression provides an upper bound on the minimum value of $d$, and that lower embedding dimensions are often sufficient in practice), which implies that we can rewrite Equation \ref{eqn:sys_conserved} as:
\begin{equation}\label{eqn:sys_conserved_delay}
    \dot{\mathbf{x}}_c=\mathbf{g}(\mathbf{x}_c,\mathbf{u},\pmb{\phi}(\mathbf{z}_e)|\pmb{\alpha},\pmb{\beta})
\end{equation}
where every term in the equation is known except for the function $\pmb{\phi}$ and parameters $\pmb{\beta}$.
Assuming that the constitutive law is a continuous function of the state variables, then we can approximate Equation \ref{eqn:sys_conserved_delay} using the following expression:
\begin{equation}\label{eqn:sys_ann}
    \dot{\mathbf{x}}_c\approx \mathbf{g}(\mathbf{x}_c,\mathbf{u},\hat{\pmb{\phi}}_{ANN}(\mathbf{z}_e)|\pmb{\alpha},\hat{\pmb{\beta}})
\end{equation}
where $\hat{\pmb{\beta}}$ is an estimate of $\pmb{\beta}$ and $\hat{\pmb{\phi}}_{ANN}$ is a feedforward artificial neural network approximation of $\pmb{\phi}$ with at least one hidden layer (in practice there will typically be several hidden layers).
Equation \ref{eqn:sys_ann} is a delay differential equation (DDE), which we will use to construct a GB-ANN.
Alternatively to the DDE formulation, one could follow the example of Krischer et al., and use an autoencoder to learn a sufficient set of data-driven intrinsic variables, and then evolve the model in time as a typical ODE\cite{Krischer1993}.

Often, there are phenomenological approximations available for the constitutive law that are functions only of the known physical variables.
For example, the Monod equation gives the growth rate of a microbe as a function of substrate concentration $S$, with two parameters, the maximum growth rate $\mu_m$ and the half saturation constant $K_s$\cite{Monod1949}.
Approximate relationships like these can be incorporated into the GB-ANN modeling framework.
We can rewrite Equation \ref{eqn:sys_ann} as follows:
\begin{equation}\label{eqn:sys_ann_prior}
    \dot{\mathbf{x}_c} \approx \mathbf{g}(\mathbf{x}_c,\mathbf{u},\mathbf{p}(\mathbf{x}_c,\mathbf{u})\hat{\pmb{\phi}}_{ANN}(\mathbf{z}_e)|\pmb{\alpha},\hat{\pmb{\beta}})
\end{equation}
where we denote the phenomenological approximation using $\mathbf{p}(\mathbf{x}_c,\mathbf{u})$.
Now, the relationship learned by the ANN model is not the complete constitutive law, but only the ``correction factor''  (in some sense, the ``overall effectiveness factor'') of how the underlying true law and the phenomenological model differ. (In other contexts, it may be useful to use $\left(\mathbf{p}(\mathbf{x}_c,\mathbf{u})+\hat{\pmb{\phi}}_{ANN}(\mathbf{z}_e)\right)$, an additive correction, instead of $\left(\mathbf{p}(\mathbf{x}_c,\mathbf{u})\hat{\pmb{\phi}}_{ANN}(\mathbf{z}_e)\right)$, depending on the nature of the prior information. We make our choice because in our application of interest---cell growth and biotechnology---it is likely that we know some, but not all factors that inhibit cell growth, and a common modeling approach is to multiply all growth inhibition factors together into one overall ``effective'' factor.)

\subsection{GB-ANN Model Structure}
\begin{figure}
    \centering
    \includegraphics[width=0.75\textwidth]{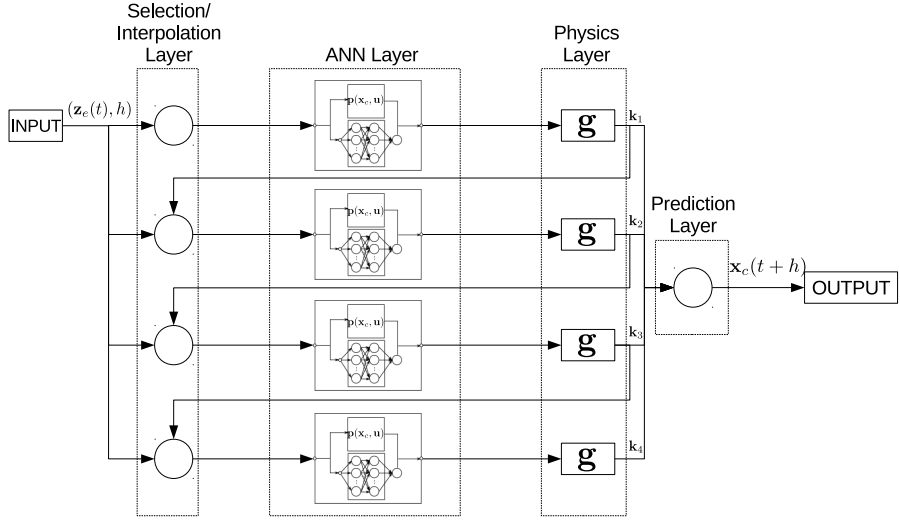}
    \caption{Schematic of a GB-ANN model. The Selection/interpolation layers sample the history of $(\mathbf{x}_c, \mathbf{u})$ at the necessary times to input to the ANN layer (see text); the ANN layer estimates the constitutive law, $\mathbf{p}(\mathbf{x}_c,\mathbf{u})\hat{\pmb{\phi}}_{ANN}(\mathbf{z}_e)$ (note that absent a phenomenological approximation, we can set $\mathbf{p}(\mathbf{x}_c,\mathbf{u})=1$), the physics layer applies the function $\mathbf{g}(\mathbf{x}_c,\mathbf{u},\mathbf{p}(\mathbf{x}_c,\mathbf{u})\hat{\pmb{\phi}}_{ANN}(\mathbf{z}_e)|\pmb{\alpha},\hat{\mathbf{\beta}})$ (time shifted as needed for the fourth order RK algorithm), the prediction layer applies the RK equation (Equation \ref{eqn:RK4}). At training time, the weights of the ANN layer and the unknown parameters, $\beta$, of the physics layer are adjustable.}
    \label{fig:gb_ann_schematic}
\end{figure}

The DDE in Equation \ref{eqn:sys_ann} is a continuous time model; however, in most scenarios data is available in discrete time intervals with sampling time of $t_s$.
Therefore, we construct a recurrent network architecture templated on a fourth order Runge-Kutta (RK) integrator that predicts the time evolution of the quantities that participate in our conservation laws\cite{Rico-Martinez1992,Rico-Martinez1993,Rico-Martinez1994}.
By approximating the DDE using the 4\textsuperscript{th} order RK approximation, and assuming that the actuators are held constant during each time step, we can write:
\begin{equation}\label{eqn:RK4}
\begin{aligned}
    \mathbf{k}_1&=\mathbf{g}\left(\mathbf{x}_c(t),\mathbf{u}(t),\mathbf{p}\left(\mathbf{x}_c(t),\mathbf{u}(t)\right)\hat{\mathbf{\phi}}_{ANN}\left(\mathbf{z}_e(t)\right)|\pmb{\alpha},\hat{\pmb{\beta}}\right)\\
    \mathbf{k}_2 &=\mathbf{g}\left(\mathbf{x}_c(t)+\frac{\mathbf{k}_1}{2},\mathbf{u}(t),\mathbf{p}\left(\hat{\mathbf{x}}_c(t+0.5h),\mathbf{u}(t)\right)\hat{\pmb{\phi}}_{ANN}\left(\hat{\mathbf{z}}_{e{\mathbf{k}_1}}(t+0.5h)\right)|\pmb{\alpha},\hat{\pmb{\beta}}\right)\\
    \mathbf{k}_3 &=\mathbf{g}\left(\mathbf{x}_c(t)+\frac{\mathbf{k}_2}{2},\mathbf{u}(t),\mathbf{p}\left(\hat{ \mathbf{x}}_c(t+0.5h),\mathbf{u}(t)\right)\hat{\pmb{\phi}}_{ANN}\left(\hat{\mathbf{z}}_{e{\mathbf{k}_2}}(t+0.5h)\right)|\pmb{\alpha},\hat{\pmb{\beta}}\right)\\
    \mathbf{k}_4 &=\mathbf{g}\left(\mathbf{x}_c(t)+\mathbf{k}_3,\mathbf{u}(t),\mathbf{p}\left(\hat{ \mathbf{x}}_c(t+h),\mathbf{u}(t)\right)\hat{\pmb{\phi}}_{ANN}\left(\hat{\mathbf{z}}_{e{\mathbf{k}_3}}(t+h)\right)|\pmb{\alpha},\hat{\pmb{\beta}}\right)\\
    \mathbf{x}_c(t+h)&=\mathbf{x}_c(t)+\frac{h}{6}\left(\mathbf{k}_1+2\mathbf{k}_2+2\mathbf{k}_3+\mathbf{k}_4\right)\\
\end{aligned}
\end{equation}
where $h$ is the step size and $\hat{\mathbf{z}}_{e{\mathbf{k}_i}}(t+h)$ is found by linear interpolation of the sampled $\{\mathbf{u}(t), \mathbf{x}_c(t)\}$ history with $\mathbf{k}_i$ used as a gradient to estimate future values of $\mathbf{x}_c$.
While our particular ANN architecture requires that $h=t_s$, we emphasize that we could (after learning the RHS of the model) create ANN architectures based on any established DDE integrator and use variable step sizes.
Equation \ref{eqn:RK4} implies that there are four passes (recurrences) through the neural network layer ($\hat{\pmb{\phi}}_{ANN}$) and the ``physics layer'' ($\mathbf{g}$) to predict $\mathbf{x}_c(t+h)$, which is illustrated schematically in Figure \ref{fig:gb_ann_schematic}.

For the ANN used to represent the unknown constitutive law, we use 3 hidden layers of 20 neurons each with the softplus activation functions \cite{glorot2011}, followed by an output layer with neurons corresponding to each element of $\mathbf{x}_c$ and linear activation functions.
We \emph{do not} use any additional regularization tools such as dropout or shrinkage penalties, as we intend to showcase how the physics layers effectively regularize the networks.
Of course, the particular choices of ANN activation function, depth, width, and regularization may be treated as tuning parameters for an individual problem.

\subsection{GB-ANN Model Training}
The GB-ANN model is trained on time series data collected at uniform time intervals $\tau$.
Trainable parameters include the weights of the ANN that represents the constitutive law as well as the parameters, $\pmb{\beta}$, used by the physics layer.
The first $d$ time steps are required to construct the delay embedding ($\mathbf{z}_e$); afterward, each new time step provides a new data point.
Given the input $\mathbf{z}_e(t)$, the GB-ANN is used to predict the quantities that participate in conservation laws at the next time step, i.e., $\mathbf{x}_c(t+\tau)$.
In principle, the GB-ANN model could be trained using any input sequence; in practice, we recommend a series of random step changes to explore the controllable subspace of the system (subject to any physical or safety constraints).
The training data are normalized (transformed to zero mean and unit variance), which requires that the physics equations be rescaled  as well.
If \emph{a priori} estimates of the mean and variance of $\pmb{\phi}$ are known, then the physics equations should be written to use normalized values of $\pmb{\phi}$ as well.

For the numerical examples below, the GB-ANN models were constructed using the Keras API with the TensorFlow backend.\cite{tensorflow2015-whitepaper, chollet2015keras}
Custom layers were designed for the selection/interpolation and prediction layers, and physics layers were developed for each of the example systems described below (see GB-ANN software package at \url{http://bitbucket.org/rlovelett/grey_box}).
Gradients were found (via automatic differentiation in TensorFlow) using back-propagation, and the models were trained using the Adam optimizer\cite{Kingma2014} using the mean squared error cost function and batches of 50 training points.

\section{Results and Discussion} \label{sec:res_disc}
We designed GB-ANNs for developing dynamical systems models of chemical/bio reactors in biotechnology, especially focusing on using light as an actuator via optogenetics.
First, as a validation step, we use a simplified isothermal model of a CSTR/chemostat that contains a single biochemical reaction subject to enzyme kinetics.
We show how to use the GB-ANN modeling framework to learn the underlying rate law, even when some of the species data are unavailable.
Second, we develop a semi-empirical mechanistic model of a bioreactor where a microbe is controlled using optogenetics.
We use this model to generate artificial training/validation data for our grey-box modeling framework.
We demonstrate how observations of only \emph{extracellular} species on which we can write balance equations (and which are measurable, even if only off-line) are sufficient for identifying the dynamical model using the GB-ANN approach, even though the underlying system also depends on \emph{intracellular} species that cannot be measured (while transcriptomic/proteomic/metabalomic measurements are possible in principle, they are impractical for dynamical system identification).

\subsection{Chemostat with a Biochemical Reaction}

For a simple demonstration of the grey-box modeling framework, we present a model of a single reaction (\ce{S ->[E] P}) in an isothermal continuous stirred tank reactor (CSTR) (i.e., a chemostat).
By writing a species balance and assuming the Michaelis-Menten rate law, the model equations are:
\begin{equation}\label{eqn:michaelis_model}
\begin{aligned}
    \dot S &=\frac{1}{\theta}(S_0-S) - \frac{k_{cat}ES}{K_M+S}\\
    \dot E &=\frac{1}{\theta}(E_0-E)\\
    P(t)&=S_0-S
\end{aligned}
\end{equation}
where the state variables $\{S,E\}$ are the substrate and enzyme concentrations, respectively, $P$ is the product concentration, the parameters $\{S_0, \theta, k_{cat}, K_M\}$ are the inlet substrate concentration, residence time, turnover number, and Michaelis-Menten constant, respectively, and $E_0$ is the inlet enzyme concentration.
Parameter values can be found in Table \ref{tab:MM_params} (note that symbols for variables used in this section are distinct from the subsequent section).

\begin{table}[ht]
    \centering
    \begin{tabular}{|c|c|c|}
        \hline
         Parameter & Symbol & Value  \\
         \hline
         Inlet Substrate Concentration & $S_0$ & \SI{3.E-2}{\mole \per \liter}\\
         Residence Time & $\theta$ & \SI{20.}{\second}\\
         Turnover Number &$k_{cat}$ & \SI{0.14}{\per \second} \\
         Michaelis-Menten Constant & $K_M$ & \SI{1.5E-2}{\mol\per\liter} \\
         \hline
    \end{tabular}
    \caption{Parameter definitions and values used in the model of the stirred tank reactor.}
    \label{tab:MM_params}
\end{table}

Now, consider that we have access to measurements of the substrate concentration, $S$, and that we can manipulate the inlet concentration of enzyme, $E_0$.
Furthermore, we will restrict our knowledge of the underlying system to the first row of Equation \ref{eqn:michaelis_model}.
Therefore, we write:
\begin{equation}\label{eqn:grey_michaelis_model}
    \dot{x}_c=\frac{1}{\alpha_1}(\alpha_2-x_c)+\phi(\mathbf{z}_e)
\end{equation}
where we have defined $x_c := S$, $u := E_0$, $\pmb{\alpha} := \{\theta, S_0\}$ (which are known constants), $\mathbf{z}_e$ the delay embedding is defined in Equation \ref{eqn:delay_embedding_definition}, and $\phi$ will be approximated using an artificial neural network. In this example, we use a delay time of \SI{10.}{\second} and embedding dimension of 5.

\begin{figure}
    \centering
    \includegraphics[width=0.75\textwidth]{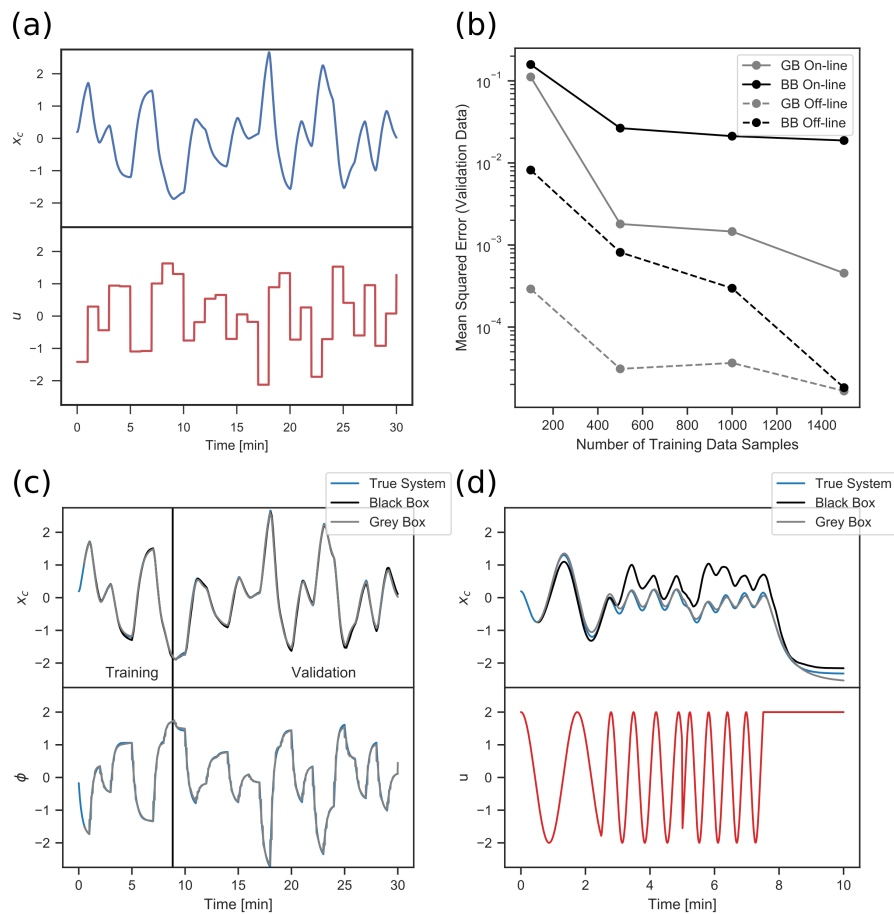}
    \caption{(a) Training and validation data. $x_c$ is substrate concentration and $u$ is inlet enzyme concentration. (b) Mean squared error of GB-ANN and BB-ANN models calculated off-line and on-line (see text) with various number of training points available; data points show the median MSE resulting from replicating the training procedure 20 times. (c) True system and modeled systems (representative example of GB-ANN and BB-ANN models) trained using 500 points, plotted with the training data and validation data series; True system and GB-ANN model system showing, $\phi$, the reaction rate. Because both models closely approximate the true system, all curves are nearly overlapping. (d) On-line performance of GB-ANN and BB-ANN for models trained with 500 points when provided an different input function than was used for training/validation, plotted with the true time series from Equation \ref{eqn:michaelis_model}. Note: all quantities are unitless as the data have been rescaled to mean zero and unit standard deviation.}
    \label{fig:MM_demo}
\end{figure}

We design a GB-ANN model as described in Section \ref{sec:theory} using Equation \ref{eqn:grey_michaelis_model} as the ``physics'' layer of the model, and do not include any phenomenological approximation for the constitutive law (i.e., $p(x_c,u)=1$ in Equation \ref{eqn:sys_ann_prior}).
Additionally, we design a simple discrete-time black-box ANN (BB-ANN) model that is given the same inputs and outputs ($\mathbf{z}_e(t)$ and $x_c(t+\tau)$, respectively) as the GB-ANN.
The BB-ANN uses three hidden layers of 20 neurons each using the softplus activation function, and an output layer with a linear activation function to predict $x_c(t+\tau)$.

By integrating Equation \ref{eqn:michaelis_model} with a series of (randomly selected) step changes in the input $S_0$, we generate a time series of 1800 points to use as training and validation data, which are plotted in Figure \ref{fig:MM_demo}a (the sampling time is 1 second).
To investigate the relative performance of the GB-ANN and BB-ANN models, we varied the number of points available for training and examined the prediction errors.
For both models, we trained the network weights for 100 epochs using the Adam algorithm with batch sizes of 50 training samples.
Figure \ref{fig:MM_demo}b shows the mean squared validation error for the grey-box and black-box models calculated using two methods: (1) ``off-line'' error is the error assuming we are given each labeled input and output, and find the mean square difference between the prediction and the label (i.e., like a typical supervised learning regression problem) and (2) ``on-line'' where the model is only provided the initial condition and the sequence of inputs and then subsequently integrated to generate a time series and the mean-squared difference between the time series generated by the model and the validation section of the data is calculated.
The time series produced by the GB-ANN and BB-ANN (trained on 500 data points) are plotted in Figure \ref{fig:MM_demo}c; we also plot, for only the GB-ANN, the output of the estimated constitutive law, $\hat \phi(z_e)$.
We observe dramatically reduced off-line error for the GB-ANN versus the BB-ANN, which is expected given that the GB-ANN incorporates previous knowledge.
The reduction in off-line error (relative to the BB-ANN) is most significant when there are comparatively fewer data points available, emphasizing how the GB-ANN model is particularly useful when producing accurate data is expensive. 
For on-line error, we see a lower mean-squared error for the GB-ANN compared with the BB-ANN model, even with up to 1500 training data points.
We also notice that the model of the constitutive law closely follows the underlying ``true'' value that is found using the Michaelis Menten law (lower panel of Figure \ref{fig:MM_demo}c).
Finally, we tested the models (GB-ANN and BB-ANN trained on 500 data points) on an entirely different input function---we assume that $E_0$ is fluctuating sinusoidally with frequency changing in time, which is shown in Figure \ref{fig:MM_demo}d.
While both models perform similarly at low frequencies (at least to a visual approximation), we observe that as the frequency increases, the GB-ANN model more closely matches the true system, whereas systematic errors remain in the BB-ANN model.

\subsection{Light-Controlled Continuous Bioreactor}\label{sec:cont_bio}
The previous example demonstrates how to design a GB-ANN for a simple system where we artificially limited our knowledge of the system dynamics by restricting our information to only one of the two balance equations.
Here, we return to the example that was introduced in Section \ref{sec:theory} where the underlying system state, and thus the relevant dynamical equations, are unknown.

We examine continuous production of a desired product such as a biofuel using a more complex chemostat model containing species growth and consumption/production reactions actuated by a light input.
This example was inspired by current research that aims to use optogenetics to control gene expression in microbes in order to, for example, increase the yield of biofuels or other products\cite{Milias-Argeitis2016,Zhao2018}.
By using optogenetics, the concentrations of different enzymes inside microbe can be influenced by applying particular wavelengths of light, which can enable dynamic control of the products that are catalyzed by these enzymes.
Using optogenetics to manipulate gene expression has resulted in record yields of the second-generation biofuel isobutanol produced in \textit{Saccharomyces cerevisiae}\cite{Zhao2018}, and has been proposed for a larger class of systems\cite{Lalwani2018}.
The example illustrates the utility of our approach when the actuator indirectly influences the quantities governed by conservation laws, in this case, by altering the unknown and unmeasured intracellular environment.

\subsubsection{Mechanistic Model of Light-Controlled Bioreactor}\label{sec:mech_bioreactor}
To test the GB-ANN modeling approach, we developed a multiscale mechanistic model to describe microbial growth and production to generate training/testing data.
The full model consists of eight ODEs, four of which are macroscopic species balances on the bioreactor, and four of which are related to intracellular concentrations of enzymes, mRNAs, and how they interact with the actuator (light input).

Let $\mathbf{x}=\{X, A, B, S, E_A, E_B, R_A, R_B\}$, where $X$ is concentration of cells,  $A$ and $B$ are extracellular concentrations of products, $S$ is concentration of a limiting nutrient, $E_A$ and $E_B$ are intracellular concentrations of enzymes that catalyze production of the products, and $R_A$ and $R_B$ are concentrations of mRNAs that code for enzymes $E_A$ and $E_B$, respectively.
The model equations are:
\begin{equation}\label{eqn:chemostat_model}
    \begin{aligned}
        \dot X &= (\mu - D)X\\
        \dot A &= \mu_{A} X - D A\\
        \dot B &= \mu_{B} X - D B\\
        \dot S &= -\left(\frac{\mu}{Y_X} + \frac{\mu_A}{Y_A} + \frac{\mu_B}{Y_B}\right)X +(S_0 -S)D\\
        \dot E_A &= -\frac{1}{\tau_1}(E_A-R_A)\\
        \dot E_B &= -\frac{1}{\tau_2}(E_B-R_B)\\
        \dot R_A &= -\frac{1}{\tau_3}(R_A-R_{A0})\\
        \dot R_B &= -\frac{1}{\tau_4}(R_B-R_{B0})
    \end{aligned}
\end{equation}
where:
\begin{equation}\label{eqn:chemostat_model_algebraic}
    \begin{aligned}
        \mu &\equiv \mu_0 S \frac{\exp{\left(-\frac{A}{K_A}-\frac{B}{K_B}-\frac{E_B}{K_{EB}}\right)}}{K_S+S}\\
        \mu_A &\equiv E_A S \frac{\exp{\left(-\frac{A}{K_{AA}}\right)}}{K_{SA}+S}\\
        \mu_B &\equiv E_B S \frac{\exp{\left(-\frac{B}{K_{BB}}\right)}}{K_{SB}+S}\\
        R_{A0} &\equiv R_{A1}+(R_{A2}-R_{A1})\frac{u^{N_A}}{K_{HA}+u^{N_A}} \\ 
        R_{B0} &\equiv R_{B1}+(R_{B2}-R_{B1})\frac{(1-u)^{N_B}}{K_{HB}+(1-u)^{N_B}}
    \end{aligned}
\end{equation}
where $u$ is the actuator (light intensity) and all parameters for Equations \ref{eqn:chemostat_model}--\ref{eqn:chemostat_model_algebraic} are provided in Table \ref{tab:chemostat_params}.
The model is based on the assumption of a Monod-like growth rate that is inhibited by products and by metabolic burden of producing enzyme B.
The light input ``switches'' the microbe from producing mRNAs corresponding to enzymes that produce species A (light on) to producing mRNAs corresponding to enzymes that produce species B (light off).
We assume first order kinetics associated with the changes to the steady-state levels of mRNAs and enzymes, and assume that the steady states change with light intensities subject to Hill functions.
Enzyme and mRNA concentrations have units of \si{\per \hour} and are understood as concentrations of enzymes (or mRNAs) that correspond to a particular specific maximum production rate (i.e., $V_{max}$, in the enzyme kinetics literature).

\begin{table}[ht]
    \centering
    \begin{tabular}{|c|c|c|}
        \hline
         Parameter & Symbol & Value  \\
         \hline
         Dilution Rate & $D$ & \SI{0.05}{\per \hour}\\
         Inlet Nutrient Concentration & $S_0$ & \SI{20.}{\gram \per \liter}\\
         Biomass Yield Coefficient & $Y_X$ & \SI{0.435}{} \\
         Species A Yield Coefficient &$Y_A$ & \SI{0.607}{} \\
         Species B Yield Coefficient &$Y_B$ & \SI{0.3}{}\\
         Nominal Specific Growth Rate & $\mu_0$ & \SI{0.22}{\per\hour}\\
         Monod Constant for Growth & $K_S$ & \SI{1.03}{\gram \per \liter}\\
         Species A  Growth Inhibition Constant & $K_A$ & \SI{7.12}{\gram \per \liter}\\
         Species B Growth Inhibition Constant & $K_B$ & \SI{0.712}{\gram \per \liter}\\
         Enzyme B Growth Inhibition Constant & $K_{EB}$ & \SI{0.5}{\per\hour}\\
         Minimum Enzyme A Production Rate & $R_{A1}$ & \SI{0}{}\\
         Maximum Enzyme A Production Rate & $R_{A2}$ & \SI{1.79}{\per\hour}\\
         Monod Constant for Production of A & $K_{SA}$ & \SI{1.68}{\gram\per\liter}\\
         Inhibition Constant for Production of A & $K_{AA}$ & \SI{14.}{\gram\per\liter}\\
         Minimum Enzyme B Production Rate & $R_{B1}$ & \SI{0.0985}{\per\hour}\\
         Maximum Enzyme B Production Rate & $R_{B2}$ & \SI{0.448}{\per\hour}\\
         Monod Constant for Production of B & $K_{SB}$ & \SI{1.68}{\gram\per\liter}\\
         Inhibition Constant for Production of B & $K_{BB}$ & \SI{14.}{\gram\per\liter}\\
         Time Constant for Enzyme A & $\tau_1$ & \SI{5}{\hour}\\
         Time Constant for Enzyme B & $\tau_2$ & \SI{6.}{\hour}\\
         Time Constant for mRNA A & $\tau_3$ & \SI{1.}{\hour}\\
         Time Constant for mRNA B & $\tau_4$ & \SI{1.}{\hour}\\
         Hill Exponent for Species A & $N_A$ & \SI{2.7}{}\\
         Hill Exponent for Species B & $N_B$ & \SI{2.3}{}\\
         Half Saturation Constant for Species A & $K_{HA}$ & \SI{0.08}{} \\
         Half Saturation Constant for Species B & $K_{HB}$ & \SI{0.30}{} \\
         \hline
    \end{tabular}
    \caption{Parameter definitions and values used in the mechanistic model of the chemostat (values were chosen such that plausible physiological behavior was observed).}
    \label{tab:chemostat_params}
\end{table}

\subsubsection{GB-ANN Model of Light-Controlled Bioreactor}
The model equations used for the GB-ANN include the four macroscopic species balances on the bioreactor (biomass $X$, chemical species $A$ and $B$, and substrate $S$), under the assumption that we do not have access to the other four states that correspond to intracellular species (enzymes $E_A$ and $E_B$, and mRNAs $R_A$ and $R_B$):
\begin{equation}\label{eqn:GB_chemostat}
    \begin{aligned}
    \dot{x}_{c1} &= \left(p_1(x_c,u)\phi_1(z_e)-\alpha_1\right)x_{c1}\\
    \dot{x}_{c2} &= -x_{c2}\alpha_1+p_2(x_c,u)\phi_2(z_e)x_{c2}\\
    \dot{x}_{c3} &= -x_{c3}\alpha_1+p_3(x_c,u)\phi_3(z_e)x_{c3}\\
    \dot{x}_{c4} &= \left(\alpha_2-x_{c4}\right)\alpha_1-\left(\frac{p_1(x_c,u)\phi_1(z_e)}{\beta_1}+\frac{p_2(x_c,u)\phi_2(z_e)}{\beta_2}+\frac{p_3(x_c,u)\phi_3(z_e)}{\beta_3}\right)x_{c1}\\
    \end{aligned}
\end{equation}
where all variables were defined in Section \ref{sec:theory}, and we emphasize again that $z_e$ is a time delay embedding.
For this example, we used a time delay of \SI{5.}{\hour} and an embedding dimension of 10.
The GB-ANN model can be compared with the mechanistic model from Section \ref{sec:mech_bioreactor}.
In terms of physical variables, $\mathbf{x}_c =\{X, A, B, S\}$, $(\mathbf{p}\pmb{\phi}) = \{\mu, \mu_A, \mu_B\}$, $\pmb{\alpha} = \{D, S_0 \}$, and $\pmb{\beta}=\{Y_X, Y_A, Y_B\}$.

We consider two versions of the model: GB1-ANN, which uses $p_1=p_2=p_3=1$ and therefore assumes no previous information about the constitutive law, and GB2-ANN, which uses the following approximate phenomenological models in conjunction with the ANN model for the constitutive law:
\begin{equation}\label{eqn:phenom_chemostat}
    \begin{aligned}
    p_1(x_c,u) &= \frac{0.18 \exp{\left(-\frac{x_{c2}}{7}-\frac{x_{c3}}{7}\right)}x_{c4}}{1+x_{c3}}\\
    p_2(x_c,u) &= \frac{0.9 \exp{\left(-\frac{x_{c2}}{10}\right)}x_{c4}}{1.5+x_{c4}}\\
    p_2(x_c,u) &= \frac{0.25  \exp{\left(-\frac{x_{c2}}{10}\right)}x_{c4}}{1.5+x_{c4}}.
    \end{aligned}
\end{equation}
Comparing Equations \ref{eqn:GB_chemostat}--\ref{eqn:phenom_chemostat} above with Equations \ref{eqn:chemostat_model}--\ref{eqn:chemostat_model_algebraic} from Section \ref{sec:mech_bioreactor} shows how these phenomenological models are missing terms, have inaccurate parameters, and incorrectly assume instantaneous dynamics, all of which require the ANN model for correction.

\begin{figure}
    \centering
    \includegraphics[width=0.65\textwidth]{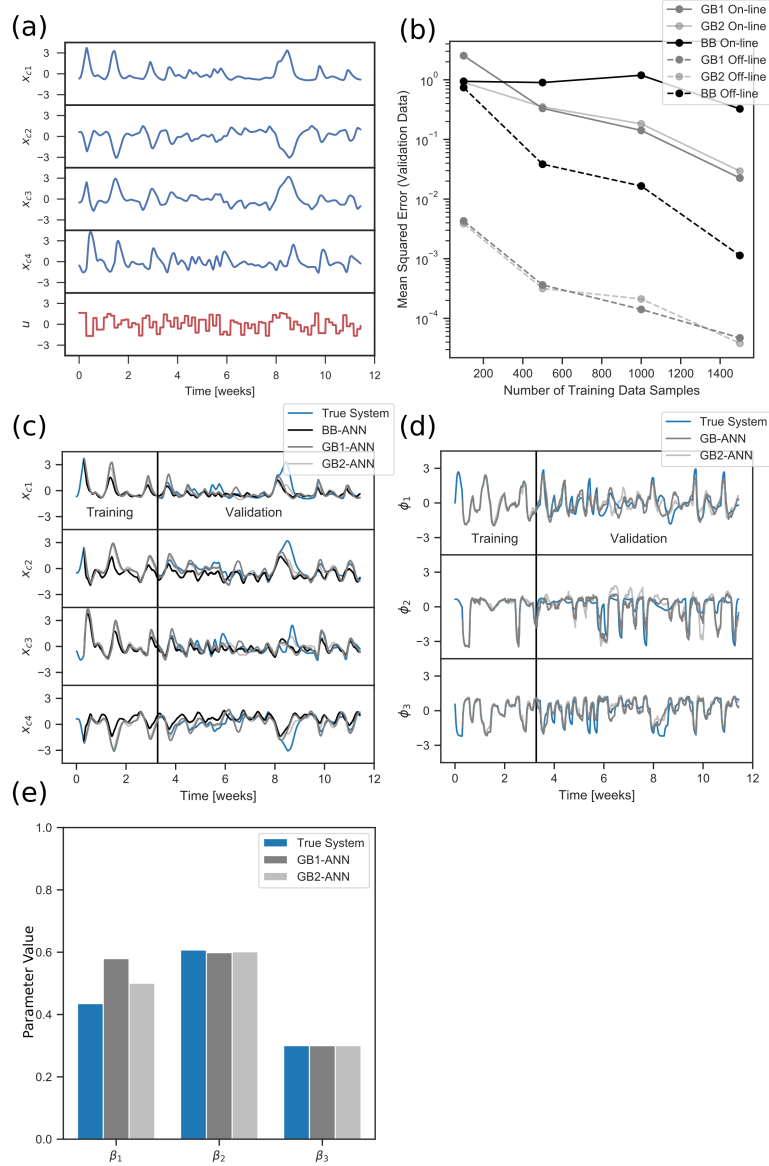}
    \caption{(a) Training and validation data. $\mathbf{x}_c$ is a vector of the species governed by conservation laws and $u$ is the actuator. (b) Mean squared error of GB1-ANN, GB2-ANN, and BB-ANN models calculated off-line and on-line (see text) with various number of training points available; data points show the median MSE resulting from replicating the training procedure 10 times. (c) True system and modeled systems (GB1-ANN, GB2-ANN, and BB-ANN models) trained using 500 points, plotted with the training data and validation data series. (d) True system and GB-ANN1 and GB-ANN2 models system showing showing the constitutive law outputs $\pmb{\phi}$. (e) The unknown constants in each of the grey-box models (i.e., yield coefficients, $\pmb{\beta}$) after training, in comparison to the true values used to generate the underlying data. Note: all time-varying quantities are unitless because the data have been rescaled to mean zero and unit standard deviation.}
    \label{fig:chemostat_demo}
\end{figure}

Now, we use these two GB-ANN models and a BB-ANN model of the system to investigate how they perform with limited training data.
We use the same ANN architecture as for the previous example, changing only the number of inputs and outputs to correspond to the larger number of species governed by conservation laws in Equation \ref{eqn:GB_chemostat}.
Training data were generated using the mechanistic model from Section \ref{sec:mech_bioreactor}, and are plotted in Figure \ref{fig:chemostat_demo}a.
The models were trained using the Adam algorithm for 1000 epochs with varying numbers of training points.

Figure \ref{fig:chemostat_demo}b shows reduced mean-squared validation errors for the GB-ANN models in comparison with the BB-ANN models for any amount of training data, except for the smallest training set, for which no models are highly accurate on-line.
GB1-ANN and GB2-ANN have comparable mean-squared errors, even with a small number of data points.
The GB-ANN models more closely track the species governed by the balance equations than the BB-ANN models (Figure \ref{fig:chemostat_demo}c), both GB-ANN models track the underlying constitutive laws  closely (Figure \ref{fig:chemostat_demo}d), and both find accurate estimates of the unknown parameters $\beta$ (Figure \ref{fig:chemostat_demo}e).

\section{Conclusions}\label{sec:conc}
In this paper, we presented an approach for physics-informed grey-box system identification that combines physics-based modeling equations with data-driven constitutive laws for applications in biotechnology.
In particular, we focus on the case where the dynamics of macroscale variables can be written using conservation laws, but they depend on constitutive laws that may be functions of microscopic/intracellular state variables.
We use an idea from nonlinear dynamics---embedology, and, in particular, time delay embeddings---to reconstruct the unknown state and write the constitutive laws as functions of known variables and their histories.
By investigating two example chemostat systems---a single biochemical reaction in a CSTR and a continuous bioreactor actuated using light---we demonstrate how incorporating the physics equations into an RK ANN effectively regularizes and reduces prediction error in the GB-ANN in comparison with a BB-ANN model with approximately the same number of trainable parameters. (Because constants $\beta$ are trainable, the GB-ANN may have slightly more trainable parameters than the BB-ANN.)

The use of time delay embeddings that are diffeomorphic to the intrinsic state of the system has consequences that could be valuable for tasks related to system identification.
As various embeddings theorems\cite{Takens1981, Whitney, Stark1997,Stark1999,Stark2003} show, different observations (with their histories) of the same system all lie on the same underlying manifold.
We can consider how these different observations and their histories (say on-line vs. off-line data; or ``cheap'' sensor data vs. ``expensive'' -omics data) could be mapped to this same manifold using (for example) manifold learning algorithms like diffusion maps\cite{Coifman2005, Coifman2006}.
We can use this knowledge to construct ``observers'' to determine the intrinsic state regardless of which information (so long as it is rich enough) is available at any given time point.
In cases of modeling uncertainty or high noise levels, we could also use nonlinear extensions to classical Kalman filters to assimilate data in real-time\cite{Tatiraju1998, Cao2004,Julier2004,Xiong2006,Guo2015}.

The GB-ANNs developed in this work are delay-differential equation (DDE) models using only known macroscopic physical quantities as state variables.
The RHS of the model equations requires only these physical quantities, actuator values, and the histories of each.
These models can therefore be deployed for on-line state estimation, model predictive control, or other process systems engineering tasks.
Because they are nonlinear, they are especially valuable for processes that may not operate at steady state, and could be useful for tasks such as economic model predictive control that maximize an economic objective rather than maintain a constant set point\cite{ellis2014tutorial}.

Using GB-ANNs constrains the model to follow physical laws, which may have benefits outside of improved training in ``medium data'' environments.
For predictive control applications, following physical laws reduces the chance that the model will results in physically impossible (and potentially dangerous) outputs.
By incorporating the known system equations into the model, it improves interpretability of the model without sacrificing the capability of deep networks to represent highly complicated functions.

The importance of reconstructing the system state from delay embeddings arises in the light-actuated bioprocesses investigated here because the system state depends on intracellular information that is unlikely to be available to the modeler in an explicit form.
The direct effect of light is to change the conformation of intracellular species, which results in a cascade of effects eventually influencing the macroscale variables.
The ubiquity of complex systems partly governed by unobserved microscale phenomena that drive macroscale variables make this approach promising across the physical and biological sciences, and we anticipate that it will find use in a number of practical domains.

It has been proposed to use optogenetics for on-line control of bioprocesses that produce useful products like biofuels using light as an actuator\cite{Lalwani2018,Zhao2018}.
Although researchers have successfully used optogenetics for control of fluorescent proteins (which can be easily measured on-line) in microbial systems\cite{Milias-argeitis2015, Milias-Argeitis2016}, control of quantities which cannot be reliably measured on-line requires a system model in order to design a state estimator.
Moving forward, we envision models like the ones developed here coupled with more easily measurable real-time data, such as optical densities or fluorescent biosensors, for accurate control and optimization of light-actuated bioprocesses.

\begin{acknowledgement}
The authors thank DARPA, an ARO MURI (IGK) and the Eric and Wendy Schmidt Transformative Technologies Fund (RJL, JLA, and IGK) for partial financial support of this research.

\end{acknowledgement}

\bibliography{main.bib}
\end{document}